\documentclass[a4paper]{article}%
\usepackage{graphicx}
\usepackage{geometry}
\usepackage{amsmath}
\usepackage{amssymb}
\usepackage{amsfonts}
\usepackage{amsthm,amscd}%
\providecommand{\U}[1]{\protect\rule{.1in}{.1in}}
\def\figurename{Figure}
\makeatletter
\renewcommand{\fnum@figure}[1]{\figurename~\thefigure.}
\makeatother
\def\tablename{Table}
\makeatletter
\renewcommand{\fnum@table}[1]{\tablename~\thetable.}
\makeatother
\def \bop {\noindent\textbf{Proof. }}
\def \eop {\hbox{}\nobreak\hfill
\vrule width 2mm height 2mm depth 0mm
\par \goodbreak \smallskip}
\newtheorem{theorem}{Theorem}[section]
\newtheorem{lemma}[theorem]{Lemma}
\newtheorem{corollary}[theorem]{Corollary}
\newtheorem{proposition}[theorem]{Proposition}
\theoremstyle{definition}
\newtheorem{definition}[theorem]{Definition}

\theoremstyle{remark}
\newtheorem{remark}[theorem]{Remark}
\numberwithin{equation}{section}

\begin{document}

\title{Existence and optimality conditions for relaxed mean-field stochastic control
problems\thanks{Partially supported by the French Algerian Cooperation
Program, Tassili 13\ MDU 887. }}
\author{{Khaled Bahlali}\thanks{ Laboratoire IMATH, Universit\'{e} du Sud-Toulon-Var,
B.P 20132, 83957 La Garde Cedex 05, France. E-mail: bahlali@univ-tln.fr}
\and {Meriem Mezerdi }\thanks{ Laboratory of Applied Mathematics, University of
Biskra, Po. Box 145, Biskra (07000), Algeria. E-mail: m\_mezerdi@yahoo.fr}
\and Brahim Mezerdi\thanks{ Laboratory of Applied Mathematics, University of
Biskra, Po. Box 145, Biskra (07000), Algeria. E-mail: bmezerdi@yahoo.fr}}
\maketitle

\begin{abstract}
We consider optimal control problems for systems governed by mean-field
stochastic differential equations, where the control enters both the drift and
the diffusion coefficient. We study the relaxed model, in which admissible
controls are measure-valued processes and the relaxed state process is driven
by an orthogonal martingale measure, whose covariance measure is the relaxed
control. This is a natural extension of the original strict control problem,
for which we prove the existence of an optimal control. Then, we derive
optimality necessary conditions for this problem, in terms of two adjoint
processes extending the known results to the case of relaxed controls.

\end{abstract}

\textbf{Key words}: Mean-field stochastic differential equation; relaxed
control; martingale measure;

adjoint process; stochastic maximum principle; variational principle.

\textbf{MSC 2010 subject classifications}, 93E20, 60H30.

\section{Introduction}

In this paper, we deal with optimal control of systems driven by mean-field
stochastic differential equations (MFSDE), \noindent where the coefficients
depend not only on the state but also on its distribution. This mean-field
equation, represents in some sense the average behavior of an infinite number
of particles, see \cite{JMW, Sn} for details. Since the earlier papers
$\cite{HMC, LasLio}$, mean-field control theory has raised a lot of interest,
motivated by applications to various fields such as game theory, mathematical
finance, communications networks, management of oil ressources. Mean-field
control problems occur in many applications, such as in a continuous-time
Markowitz's mean--variance portfolio selection model where the variance term
involves a quadratic function of the expectation. The inclusion of this
mean-field terms in the coefficients introduces time inconsistency, leading to
the failure of Bellmann principle. For this kind of problems, the stochastic
maximum principle, provides a powerful tool to solve them, see \cite{BucDjeLi,
Bens} and the references therein. The first objective of the present paper is
to investigate the problem of existence of an optimal control. It is well
known that in the absence of convexity assumptions, this problem has no
optimal solution. Therefore it is natural to embedd the set of strict controls
into a wider class of measure valued controls, enjoying good compactness
properties, called relaxed controls. We show that the right state process
associated with a relaxed control, satisfies a MFSDE driven by an orthogonal
martingale measure rather than a Brownian motion.\ For this model, we prove
that the strict and relaxed control problems have the same value function and
that an optimal relaxed control exists. Our result extends in particular
\cite{BahDjeMez, EKNJ, Fl, MezBah} to mean field controls and \cite{BMM} to
the case of a MFSDE with a controlled diffusion coefficient. The proof is
based on tightness properties of the underlying processes and the Skorokhod
selection theorem. In a second step, we establish necessary conditions for
optimality in the form of a relaxed stochastic maximum principle, obtained via
the first and second order adjoint processes. This result generalizes Peng's
stochastic maximum principle \cite{Peng}, to mean field control problems and
\cite{BucDjeLi} to relaxed controls. The other advantage is that our maximum
pinciple applies to a natural class of controls, which is the closure of the
class of strict controls, for which we have existence of an optimal control.
The proof of the main result is based on the approximation of the relaxed
control problem by a sequence of strict control problems. Then Ekeland's
variational principle is applied to get necessary conditions of
near-optimality, for the sequence of near optimal strict controls. The result
is obtained by a passage to the limit in the state equation as well as in the
adjoint processes. The resulting first and second order adjoint processes are
solutions of linear BSDEs driven by a Brownian motion and an orthogonal square
integrable martingale. Moreover, our result is given via an approximation
procedure, so that it could be convenient for numerical computation.

\section{Assumptions and preliminaries}

Let $(\Omega,\mathcal{F},P)$ be a probability space$,$ equipped with a
filtration $\left(  \mathcal{F}_{t}\right)  ,$ satisfying the usual conditions
and $\left(  W_{t}\right)  $ a $\left(  \mathcal{F}_{t},P\right)  -$Brownian
motion. Let $\mathbb{A}$ be some compact metric space called the action space.
A strict control $\left(  u_{t}\right)  $ is a measurable, $\mathcal{F}_{t}-$
adapted process with values in the action space $\mathbb{A}$. We denote
$\mathcal{U}_{ad}$ the space of strict controls.

\noindent The state process corresponding to a strict control is the unique
solution, of the mean-field stochastic differential equations (MFSDE)%

\begin{equation}
dX_{t}=b(t,X_{t},E(X_{t}),u_{t})dt+\sigma(t,X_{t},E(X_{t}),u_{t})dW_{t};\text{
}X_{0}=x \label{MFSDE}%
\end{equation}

\noindent and the corresponding cost functional is given by

\begin{center}
$J(u)=E\left(  \int_{0}^{T}h(t,X_{t},E(X_{t}),u_{t}\right)  dt+g(X_{T}%
,E(X_{T})).$
\end{center}

\noindent The coefficients of the state equation as well as of the cost
functional are of mean-field type, in the sense that they depend not only on
the state process, but also on its marginal law, through its expectation.

\noindent The objective is to minimize $J(u)$ over the space $\mathcal{U}%
_{ad}$ , that is to find $u^{\ast}\in$ $\mathcal{U}_{ad}$ such that
$J(u^{\ast})=\inf\left\{  J(u),u\in\mathcal{U}_{ad}\right\}  .$

\noindent Let us consider the following assumptions which will be used in
different combinations throughout the paper.

$\mathbf{(H}_{\mathbf{1}}\mathbf{)}$ $b:\left[  0,T\right]  \times
\mathbb{R}\times\mathbb{R}\times\mathbb{A}\longrightarrow\mathbb{R}$,
$\sigma:\left[  0,T\right]  \times\mathbb{R}\times\mathbb{R}\times
\mathbb{A}\longrightarrow\mathbb{R}$ are bounded continuous functions such
that $b(t,.,.,a)$ and $\sigma(t,.,.,a)$ are Lipschitz continuous in $(x,y).$\ 

$\mathbf{(H}_{\mathbf{2}}\mathbf{)}$ $h:\left[  0,T\right]  \times
\mathbb{R}\times\mathbb{R}\times\mathbb{A}\longrightarrow\mathbb{R}$ and
$g:\mathbb{R}\times\mathbb{R}\longrightarrow\mathbb{R}$, are bounded
continuous functions such that $h(t,.,.,a)$ and $g(.,.)$ are are Lipschitz
continuous in $(x,y).$

$\mathbf{(H}_{\mathbf{3}}\mathbf{)}$ $b(t,.,.,a)$, $\sigma(t,.,.,a),$
$h(t,.,.,a)$ and $g(.,.)$ are twice continuously differentiable with respect
to $(x,y),$ and their derivatives are bounded and continuous in $(x,y,a)$.

\noindent Without loss of generality, the coefficients are assumed to be one
dimensional as in \cite{BucDjeLi}, to avoid heavy notations in the definition
of adjoint processes.

\noindent Under assumption $\mathbf{(H}_{\mathbf{1}}\mathbf{),}$ according to
\cite{JMW} Prop.1.2, for each $u\in\mathcal{U}_{ad}$ the MFSDE(\ref{MFSDE})
has a unique strong solution, such that for every $p>0$ we have $E(\left\vert
X_{t}\right\vert ^{p})<+\infty.$ Moreover the cost functional is well defined.

\section{The relaxed control problem}

\subsection{The space of relaxed controls}

As it is well known in control theory, in the absence of convexity conditions,
an optimal control may fail to exist in the set $\mathcal{U}_{ad}$ of strict
controls (see e.g. \cite{Fl}). This suggests that the set\ of strict controls
is too narrow and should be embedded into a wider class of relaxed controls,
with nice compactness properties. For the relaxed model, to be a true
extension of the original control problem, the following both conditions must
be satisfied:

\noindent i) The value functions of the original and the relaxed control
problems must be equal.

\noindent ii) The relaxed control problem must have an optimal solution.

\noindent The idea of relaxed control is to replace the $\mathbb{A}$-valued
process $(u_{t})$ with a $\mathbb{P}(\mathbb{A})$-valued process $(\mu_{t})$,
where $\mathbb{P}(\mathbb{A})$ is the space of probability measures equipped
with the topology of weak convergence. Then $(\mu_{t})$ may be identified as a
random product measure on $[0,T]\times\mathbb{A}$, whose projection on $[0,T]$
coincides with Lebesgue measure. Let $\mathbb{V}\ $be the set of product
measures $\mu$ on $[0,T]\times\mathbb{A}$ whose projection on $\left[
0,T\right]  $ coincides with the Lebesgue measure $dt$. It is clear that every
$\mu$ in $\mathbb{V}$ may be disintegrated as $\mu=dt.\mu_{t}(da)$, where
$\mu_{t}(da)$ is a transition probability. The elements of $\mathbb{V}$ are
called Young measures in deterministic theory. $\mathbb{V}$ as a closed
subspace of the space of positive Radon measures $\mathbb{M}_{+}%
([0,T]\times\mathbb{A})$ is compact for the topology of weak convergence. In
fact it can be proved that it is compact also for the topology of stable
convergence, where test functions are measurable, bounded functions $f(t,a)$
continous in $a,$ see \cite{EKNJ} for further details.

\begin{definition}
A relaxed control on the filtered probability space $\left(  \Omega
,\mathcal{F},\mathcal{F}_{t},P\right)  $ is a random variable $\mu=dt.\mu
_{t}(da)$ with values in $\mathbb{V}$, such that $\mu_{t}(da)$ is
progressively measurable with respect to $(\mathcal{F}_{t})$ and such that for
each $t$, $1_{(0,t]}.\mu$ is $\mathcal{F}_{t}-$measurable.
\end{definition}

\begin{remark}
The set $\mathcal{U}_{ad}$ of strict controls is embedded into the set of
relaxed controls by identifying $u_{t}$ with $dt\delta_{u_{t}}(da).$
\end{remark}

\subsection{The relaxed state equation}

The question now is to define the natural state process associated to a
relaxed control. In deterministic control or in the stochastic theory where
only the drift is controlled, one has just to replace in equation
(\ref{MFSDE}) the drift by the same drift integrated against the relaxed
control. Now we are in a situation where both the drift and the diffusion
coefficient are controlled. Following \cite{JMW} Prop. 1.10, the existence of
a weak solution of equation (\ref{MFSDE}) associated with a strict control $u$
is equivalent to the existence of a solution for the non linear martingale problem:

\begin{center}
$f(X_{t})-f(X_{0})-\int_{0}^{t}L^{P_{X_{s}}}f(s,X_{s},u_{s})\,ds$ \textit{is
a} $P-$\textit{martingale}$\mathit{,}$
\end{center}

\noindent\textit{for every }$f\in C_{b}^{2},$\textit{\ for each }%
$t>0,$\textit{\ where }$L$\textit{\ is the infinitesimal generator }associated
with equation (\ref{MFSDE}),

\begin{center}
$L^{\nu}f(t,x,a)=\frac{1}{2}\sigma^{2}\frac{\partial^{2}f}{\partial x^{2}%
}(t,x,a)+b\frac{\partial f}{\partial x}(t,x,a),$
\end{center}

$b=b(t,x,\left\langle y,\nu\right\rangle ,a)$ and $\sigma^{2}=\sigma
^{2}(t,x,\left\langle y,\nu\right\rangle ,a)$ where $\nu\in\mathbb{P}%
_{1}(\mathbb{R}),$ the space of probability measures on $\mathbb{R}.$

\noindent Therefore, the natural relaxed martingale problem associated to a
relaxed control is defined as follows:

\begin{center}
$f(X_{t})-f(X_{0})-\int_{0}^{t}\int_{A}L^{P_{X_{s}}}f(s,X_{s},a)\,\mu
_{s}(da)ds$ is\ a $P$-martingale \textit{for each }$f\in C_{b}^{2}%
,$\textit{\ for each }$t>0.$
\end{center}

\noindent The following theorem gives a pathwise representation of the
solution of the relaxed martingale problem, in terms of a mean-field
stochastic differential equation driven by an orthogonal martingale measure.

\begin{theorem}
\ \textit{1)Let }$P$\textit{\ be a solution of the relaxed martingale problem
. Then }$P$\textit{\ is the law of an adapted, continuous process }%
$X$\textit{\ defined on an extension of the space }$\left(  \Omega
,\mathcal{F},\mathcal{F}_{t},P\right)  $\textit{, which is a solution of the
following MFSDE:}%
\begin{equation}
dX_{t}=%
{\textstyle\int_{\mathbb{A}}}
b(t,X_{t},E\left(  \mathit{\mathit{X}_{t}}\right)  ,a)\,\mu_{t}(da)dt+%
{\textstyle\int_{\mathbb{A}}}
\sigma(t,X_{t},E(X_{t}),a)\,M(da,dt);\text{ }\mathit{X}_{0}=x \label{RMFSDE}%
\end{equation}
\textit{where }$M\ $is\textit{ an orthogonal continuous martingale measure,
with intensity }$dt\mu_{t}(da).$

\textit{2) If the coefficients }$b$\textit{\ and }$\sigma$\textit{\ are
Lipschitz in }$x$\textit{, }$\mathit{y}$\textit{, uniformly in }%
$t$\textit{\ and }$a$\textit{, equation (\ref{RMFSDE})\ has a unique pathwise
solution.}
\end{theorem}

\bop1)The proof is based essentially on the same arguments as in \cite{EM},
Theorem IV-2 and \cite{JMW}, Prop. 1.10.

2) Since the coefficients are Lipschitz continuous, then following the same
steps as in \cite{JMW} and \cite{EM}, it is not difficult to prove that
Equation (\ref{RMFSDE}) has a unique solution such that for every $p>0$ we
have $E(\sup_{t\in\left[  0,T\right]  }\left\vert X_{t}\right\vert
^{p})<+\infty.$ \eop

\begin{remark}
i) Note that the othogonal martingale measure corresponding to the relaxed
control $\mathit{dt}\mu_{t}(da)$ is not unique.

ii) From now on, the probability space is an extension of the initial
probability space. The Brownian motion $\left(  W_{t}\right)  $ remains a
Brownian motion on this new probability space, but the filtration is no longer
the natural filtration of $\left(  W_{t}\right)  .$
\end{remark}

\begin{definition}
Let $(\Omega,\mathcal{F},\mathcal{F}_{t},P)$ be a filtered probability space
and $M(t,B)$ a random process, where $B\in\mathcal{B}\left(  \mathbb{A}%
\right)  .$ $M$ is a ($\mathcal{F}_{t},P)-$martingale measure if:

1)For every $B\in\mathcal{B}\left(  \mathbb{A}\right)  ,\left(  M(t,B)\right)
_{t\geq0}$ is a square integrable martingale, $M(0,B)=0$.

2)For every $t>0$, $M(t,.)$ is a $\sigma-$finite $L^{2}$-valued measure.

It is called continuous if for each $B\in\mathcal{B}\left(  \mathbb{A}\right)
,$ $M(t,B)$ is continuous and orthogonal if $M(t,B).M(t,C)$ is a martingale
whenever $B\cap C=\phi.$
\end{definition}

\begin{remark}
When the martingale measure $M$ is orthogonal, one can prove \cite{EM} the
existence of a random positive $\sigma$-finite measure $\mu\left(
dt,da\right)  \ $on $\left[  0,T\right]  \times\mathbb{A},$ such that
$\left\langle M(.,B),M(.,B)\right\rangle _{t}=\mu\left(  \left[  0,t\right]
\times B\right)  $ for all $t>0$ and $B\in\mathcal{B}\left(  \mathbb{A}%
\right)  .$ $\mu\left(  dt,da\right)  $ is called the covariance measure of
$M$.
\end{remark}

\textbf{Example \thinspace}Let ${\mathbb{A}}=\left\{  a_{1},a_{2},\cdots
,a_{n}\right\}  $ be a finite set. Then every relaxed control $dt\,\mu
_{t}(da)$ will be a convex combination of the Dirac measures $dt\,\mu
_{t}(da)=\sum_{i=1}^{n}\alpha_{t}^{i}\,dt\,\delta_{a_{i}}(da),$ where for
every $i$, $\alpha_{t}^{i}$ is a real--valued process such that $0\leq
\alpha_{t}^{i}\,\leq1$ and $\sum_{i=1}^{n}\alpha_{t}^{i}\,=1$. It is obvious
that the solution of the relaxed martingale problem is the law of the solution
of the SDE

\begin{center}
$dX_{t}=\sum\limits_{i=1}^{n}b(t,X_{t},E\left(  X_{t}\right)  ,a_{i}%
)\alpha_{t}^{i}dt+\sum\limits_{i=1}^{n}\sigma(t,X_{t},E\left(  X_{t}\right)
,a_{i})\left(  \alpha_{t}^{i}\right)  ^{1/2}\,dW_{s}^{i},\quad\quad X_{0}=x,$
\end{center}

\noindent where the $W^{i}$ are independant Brownian motions, defined on an
extension of the initial probability space. The process $M$ given by
$M(\left[  0,t\right]  \times A)=\sum\limits_{i=1}^{n}\int_{0}^{t}\left(
\alpha_{s}^{i}\right)  ^{1/2}1_{\left\{  a_{i}\in A\right\}  }dW_{s}^{i}$ is
in fact an orthogonal continuous martingale measure (cf. \cite{EM}) with
intensity $\mu_{t}(da)dt=\sum\alpha_{t}^{i}\,\delta_{a_{i}}(da)dt$. Thus, the
last SDE can be expressed in terms of $M$ and $\mu$ as follows:

\begin{center}
$dX_{t}=\int\nolimits_{\mathbb{A}}b(t,X_{t},E\left(  X_{t}\right)
,a)\,\mu_{t}(da)dt+\int\nolimits_{\mathbb{A}}\sigma(t,X_{t},E(X_{t}%
),a)\,M(da,dt)$
\end{center}

\subsection{Approximation of the relaxed model}

The relaxed control problem is now driven by equation%

\begin{equation}
dX_{t}=%
{\textstyle\int_{\mathbb{A}}}
b(t,X_{t},E(X_{t}),a)\,\mu_{t}(da)dt+%
{\textstyle\int_{\mathbb{A}}}
\sigma(t,X_{t},E(X_{t}),a)\,M(da,dt),\text{ }X_{0}=x, \label{RMFSDE2}%
\end{equation}

\noindent and accordingly the relaxed cost functional is given by%

\begin{equation}
J(\mu)=E\left(
{\textstyle\int_{0}^{T}}
{\textstyle\int_{\mathbb{A}}}
h(t,X_{t},E(X_{t}),a)\mu_{t}(da)dt+g(X_{T},E(X_{T}))\right)  . \label{RCOST}%
\end{equation}

\noindent We show in this section that the strict and the relaxed control
problems have the same value function. This is based on the chattering lemma
and the stability of the state process with respect to the control variable.

\begin{lemma}
\textbf{(}Chattering lemma\textbf{)\ }i)\textbf{ }\textit{Let }$(\mu_{t}%
)$\textit{\ be a relaxed control}$.$\textit{\ Then there exists a sequence of
adapted processes }$(u_{t}^{n})$\textit{\ with values in }$\mathbb{A}%
$\textit{, such that the sequence of random measures }$\left(  \delta
_{u_{t}^{n}}(da)\,dt\right)  $\textit{\ converges in }$\mathbb{V}$\textit{ to
}$\mu_{t}(da)\,dt,$\textit{ }$P-a.s.$

ii) For any $g$ continuous in $\left[  0,T\right]  \times\mathbb{M}%
_{1}(\mathbb{A})$ such that $g(t,.)$ is linear, we have

$\underset{n\rightarrow+\infty}{\lim}\int_{0}^{t}g(s,\delta_{u_{s}^{n}%
})ds=\int_{0}^{t}g(s,\mu_{s})ds$ uniformly in $t\in\left[  0,T\right]  ,$
$P-a.s.$
\end{lemma}

\bop See \cite{EKNJ} and \cite{Fl} Lemma 1 page 152.\eop

\begin{proposition}
1) Let $\mu=\mu_{t}(da)\,dt$ a relaxed control. Then there exist a continuous
orthogonal martingale measure $M(dt,da)$ whose covariance measure is given by
$\mu_{t}(da)\,dt.$

2) If we denote $M^{n}(t,B)=\int\nolimits_{0}^{t}\int\nolimits_{B}%
\delta_{u_{s}^{n}}(da)dW_{s},$ where $\left(  u^{n}\right)  $ is defined as in
the last Lemma, then for every bounded predictable process $\varphi
:\Omega\times\left[  0,T\right]  \times\mathbb{A}\rightarrow\mathbb{R}$, such
that $\varphi(\omega,t,.)$ is continuous$,$ we have

$E\left[  \left(  \int\nolimits_{0}^{t}\int\nolimits_{\mathbb{A}}%
\varphi(\omega,t,a)M^{n}(dt,da)-\int\nolimits_{0}^{t}\int\nolimits_{\mathbb{A}%
}\varphi(\omega,t,a)M(dt,da)\right)  ^{2}\right]  \rightarrow0$ as
$n\longrightarrow+\infty,$

for a suitable Brownian motion $W$ defined on an eventual extension of the
probability space.
\end{proposition}

\bop See \cite{Mel} pages 196-197.\eop

\begin{proposition}
\textit{1) Let }$X_{t},$ $X_{t}^{n}$ \textit{be the solutions of state
equation (}\ref{RMFSDE}) corresponding to $\mu$ and $u^{n},$ where $\mu$ and
$u^{n}$ are defined as in the chattering lemma\textit{. Then \ }%
$\underset{n\rightarrow\infty}{\lim}E(\sup_{0\leq t\leq T}\left\vert X_{t}%
^{n}-X_{t}\right\vert ^{2})=0.$

2) Let $J(u^{n})$ and $J(\mu)$ be the expected costs corresponding
respectively to $u^{n}$ and $\mu,$ then $\left(  J\left(  u^{n}\right)
\right)  $ converges to $J\left(  \mu\right)  .$
\end{proposition}

\bop Similar to \cite{BMM}, Proposition 2. \eop

\begin{remark}
As a consequence of the last proposition, it holds that the infimum among
relaxed controls is equal to the infimum among strict controls, that is the
value functions for the relaxed and strict models are the same.
\end{remark}

\subsection{Existence of an optimal relaxed control}

The following theorem, which is the main result of this section, extends
\cite{BahDjeMez, EKNJ, Fl} to systems driven by mean field SDEs and \cite{BMM}
to mean field SDEs with controlled diffusion coefficient.

\begin{theorem}
Under assumptions $\mathbf{(H_{1})}$, $\mathbf{(H_{2})}$, there exist an
optimal relaxed control.
\end{theorem}

\bop Let $\left(  \mu^{n}\right)  _{n\geq0}$ be a minimizing sequence, that is
$\lim_{n\rightarrow\infty}J\left(  \mu^{n}\right)  =\inf_{q\in\mathcal{R}%
}J\left(  \mu\right)  $ and let $X^{n}$ be the unique solution of
(\ref{RMFSDE2}), associated with $\mu^{n}$ and $M^{n}$ where $M^{n}$ is a
continuous orthogonal martingale measure with intensity $\mu^{n}.$ We will
prove that the sequence $(\mu^{n},M^{n},X^{n})$ is tight and then show that we
can extract a subsequence, which converges in law to a process $(\widehat{\mu
},\widehat{M},\widehat{X}),$ satisfying the same MFSDE. To finish the proof we
show that the sequence of cost functionals $(J(\mu^{n}))_{n}$ converges to
$J(\widehat{\mu})$ which is equal to $\underset{\mu\in\mathcal{R}}{\inf
}J\left(  \mu\right)  $ and then we conclude that $(\widehat{\mu},\widehat
{M},\widehat{X})$ is optimal.

\textbf{Step 1}: $(\mu^{n})_{n}$ is relatively compact in $\mathbb{V}$.

\noindent The relaxed controls $\mu^{n}$ are random variables with values in
the space $\mathbb{V}$ which is compact. Then Prohorov's theorem yields that
the family of distributions associated to $(\mu^{n})_{n\geq0}$ is tight, then
$(\mu^{n})_{n}$ is relatively compact in $\mathbb{V}$.

\textbf{Step 2:} $\left(  M^{n}\right)  $ is tight in $C\left(  \left[
0,1\right]  ,\mathcal{S}^{\prime}\right)  $

\noindent The martingale measures $M^{n},$ $n\geq0,$ can be considered as
random variables with values in $C\left(  \left[  0,1\right]  ,\mathcal{S}%
^{\prime}\right)  ,$ where $\mathcal{S}$ is the Schwartz space of rapidly
decreasing functions. By \cite{Mit}, Theorem 5.1, it is sufficient to show
that for every $\varphi$ in $\mathcal{S}$ the family $\left(  M^{n}\left(
\varphi\right)  ,n\geq0\right)  $ is tight in $C\left(  \left[  0,T\right]
,\mathbb{R}^{d}\right)  $ where $M^{n}\left(  \omega,t,\varphi\right)
=\int_{\mathbb{A}}\varphi(a)M^{n}(\omega,t,da).$ Let $p>1$ and $s<t$, by the
Burkholder-Davis-Gundy inequality we have
\[
E\left(  \left\vert M_{t}^{n}(\varphi)-M_{s}^{n}(\varphi)\right\vert
^{2p}\right)  \leq C_{p}\sup\nolimits_{a\in\mathbb{A}}\left\vert
\varphi(a)\right\vert ^{2p}\left\vert t-s\right\vert ^{p}=K_{p}\,\left\vert
t-s\right\vert ^{p},
\]
where $K_{p}$ is a constant depending on $p$ and $\varphi$. Then, the
Kolmogorov tightness criteria in $C\left(  \left[  0,T\right]  ,\mathbb{R}%
^{d}\right)  $ is fulfilled and the sequence $\left(  M^{n}\left(
\varphi\right)  \right)  $ is tight. Therefore $\left(  M^{n}\right)  $ is
tight in $C\left(  \left[  0,1\right]  ,\mathcal{S}^{\prime}\right)  .$

\textbf{Step 3: }$\left(  X^{n}\right)  _{n\geq0}$\textit{\ is tight in
}$C\left(  \left[  0,T\right]  ,\mathbb{R}^{d}\right)  $

\noindent Let $p>1$ and $s<t$. Using usual arguments from stochastic calculus
and the boundness of the coefficients $b$ and $\sigma,$ it is easy to show
that
\[
E\left(  \left\vert X_{t}^{n}-X_{s}^{n}\right\vert ^{2p}\right)  \leq
C_{p}\,\left\vert t-s\right\vert ^{p}%
\]
which yields the tightness of $\left(  X_{t}^{n},n\geq0\right)  $ \textit{in
}$C\left(  \left[  0,T\right]  ,\mathbb{R}^{d}\right)  .$

\textbf{Step 4:} The sequence of processes $(\mu^{n},M^{n},X^{n})$ is tight on
the space $\Gamma=\mathbb{V\times}C\left(  \left[  0,1\right]  ,\mathcal{S}%
^{\prime}\right)  \times C\left(  \left[  0,T\right]  ,\mathbb{R}^{d}\right)
$, then by the Skorokhod representation theorem, there exists a probability
space $\left(  \widehat{\Omega},\widehat{\mathcal{F}},\widehat{P}\right)  $, a
sequence $\widehat{\gamma}^{n}=\left(  \widehat{\mu}^{n},\widehat{M^{n}%
},\widehat{X^{n}}\right)  $ and $\widehat{\gamma}=\left(  \widehat{\mu
},\widehat{M},\widehat{X}\right)  $ defined on this space such that:

(i) for each $n\in\mathbb{N}$, law$\left(  \gamma^{n}\right)  =$ law$\left(
\widehat{\gamma}^{n}\right)  $,

(ii) there exists a subsequence $\left(  \widehat{\gamma}^{n_{k}}\right)
$\ of $\left(  \widehat{\gamma}^{n}\right)  $, still denoted by $\left(
\widehat{\gamma}^{n}\right)  $, which converges to $\widehat{\gamma}%
,\widehat{P}$-$a.s$. on the space $\Gamma.$

This means in particular that the sequence of relaxed controls $(\widehat{\mu
}^{n})$ converges in the weak topology to $\widehat{\mu},$ $\widehat{P}-a.s.$
and $\left(  \widehat{M^{n}},\widehat{X^{n}}\right)  $ converges to $\left(
\widehat{M},\widehat{X}\right)  ,$ $\widehat{\mathbb{P}}-a.s.$ in $C\left(
\left[  0,1\right]  ,\mathcal{S}^{\prime}\right)  \times C\left(  \left[
0,T\right]  ,\mathbb{R}^{d}\right)  $.

\noindent According to property (i), we get

\begin{center}
$\widehat{X_{t}^{n}}=x+\int\nolimits_{0}^{t}\int\nolimits_{\mathbb{A}}b\left(
s,\widehat{X_{s}^{n}},E(\widehat{X_{s}^{n}}),a\right)  \widehat{\mu}_{s}%
^{n}(da)ds+\int\nolimits_{0}^{t}\int\nolimits_{\mathbb{A}}\sigma\left(
s,\widehat{X_{s}^{n}},E(\widehat{X_{s}^{n}}),a\right)  \widehat{M^{n}%
}(ds,da),$ $\widehat{X_{0}^{n}}=x.$
\end{center}

\noindent Since the coefficients $b,$ $\sigma$ are Lipschitz continuous in
$(x,y),$ then according to property (ii) and using similar arguments as in
\cite{Sk} page 32, it holds that the first and second terms converge in
probability to the corresponding terms without the superscript $n$. Now, since
$b$ and $\sigma$ are Lipschitz continuous, then $\widehat{X}$ is the unique
solution of the MFSDE

\begin{center}
$\widehat{X_{t}}=x+\int\nolimits_{0}^{t}\int\nolimits_{\mathbb{A}}b\left(
s,\widehat{X_{s}},E(\widehat{X_{s}}),a\right)  \widehat{\mu}_{s}%
(da)ds+\int\nolimits_{0}^{t}\int\nolimits_{\mathbb{A}}\sigma\left(
s,\widehat{X_{s}},E(\widehat{X_{s}}),a\right)  \widehat{M}(ds,da),$
$\widehat{X_{0}}=x.$
\end{center}

\noindent To finish the proof of Theorem 3.11, it remains to check that
$\widehat{\mu}$ is an optimal control.

\noindent The functions $b$ and $\sigma$ are Lipschitz continuous, then
according to the above properties (i)-(ii) we get
\begin{align*}
\underset{\mu\in\mathcal{R}}{\inf}J\left(  \mu\right)   &  =\underset
{n\rightarrow\infty}{\lim}E\left[  \int\nolimits_{0}^{T}\int
\nolimits_{\mathbb{A}}h(t,X_{t}^{n},E(X_{t}^{n}),a)\mu_{t}^{n}(da)dt+g(X_{T}%
^{n},E(X_{T}^{n}))\right] \\
&  =\underset{n\rightarrow\infty}{\lim}\widehat{E}\left[  \int\nolimits_{0}%
^{T}\int\nolimits_{\mathbb{A}}h\left(  t,\widehat{X_{t}^{n}},E(\widehat
{X_{t}^{n}}),a\right)  \widehat{\mu}_{t}^{n}(da)dt+g(\widehat{X_{T}^{n}%
},E(\widehat{X_{T}^{n}}))\right] \\
&  =\widehat{E}\left[  \int\nolimits_{0}^{T}\int\nolimits_{\mathbb{A}}h\left(
t,\widehat{X_{t}},E(\widehat{X_{t}}),a\right)  \widehat{\mu}_{t}%
(da)dt+g(\widehat{X_{T}},E(\widehat{X_{T}}))\right]  .
\end{align*}

\noindent Hence $\widehat{\mu}$ is an optimal control. \eop

\begin{remark}
The proof of the last Theorem is based on tightness and weak convergence
techniques. Then it is possible to prove it by using the non linear martingale
problem, following the same steps as in \cite{EKNJ}, without using the
pathwise representation of the solution.
\end{remark}

\section{The relaxed maximum principle}

\noindent In this section, we shall derive necessary conditions for
optimality, satisfied by an optimal relaxed control. To achieve this goal, we
begin by deriving necessary conditions for near optimality, satisfied by a
minimizing sequence of strict controls, which converges to the relaxed
control. Then we pass to the limit in the state equation as well as in the
adjoint processes to obtain the relaxed maximum principle.

\noindent Throughout this section assumptions $\mathbf{(H}_{\mathbf{1}%
}\mathbf{),(H}_{\mathbf{2}}\mathbf{)}$ and $\mathbf{(H}_{\mathbf{3}}%
\mathbf{)}$ will be in force.

\subsection{Necessary conditions for near optimality}

\noindent Let $\mu=dt\mu_{t}(da)$ be an optimal relaxed control and $X$ be the
corresponding state process, solution of \ref{RMFSDE2}. According to the
optimality of $\mu$ and the chattering lemma, there exists a sequence $\left(
u_{n}\right)  \subset\mathcal{U}_{ad}$ such that:

\noindent$J(u^{n})=J(\mu^{n})\leq\inf\left\{  J(\mu);\text{ }\mu\in
\mathcal{R}\right\}  +\varepsilon_{n},$ where $\mu^{n}=dt\delta_{u_{t}^{n}%
}\left(  da\right)  $ and $\underset{n\rightarrow+\infty}{\lim}$
$\varepsilon_{n}=0.$

\noindent In this section, we apply Ekeland's variational principle
\cite{Ekel}, to establish necessary conditions for near-optimality satisfied
by the minimizing sequence $(u^{n})$.

\begin{lemma}
\textbf{(}Ekeland\textbf{) }\textit{Let }$(V,d)$\textit{\ be a complete metric
space and }$F:V\longrightarrow\mathbb{R}\mathbf{\cup\{+}\infty\mathbf{\}}$
\textit{be lower-semicontinuous and bounded from below. Given }$\epsilon
>0,$\textit{\ suppose }$u^{\epsilon}\in V$\textit{\ satisfies }$F\left(
u^{\epsilon}\right)  \leq\inf\left\{  F\left(  v\right)  \text{ };v\in
V\right\}  +\epsilon$\textit{. Then for any }$\lambda>0$\textit{, there exists
}$v\in V$\textit{\ such that }$F\left(  v\right)  \leq F\left(  u^{\epsilon
}\right)  ,$ $d(u^{\epsilon},v)\leq\lambda$ and $\forall w\neq v$%
\textit{\ }$;$ $F(v)<F(w)+\varepsilon/\lambda.d(w,v).$
\end{lemma}

\noindent For $u,$ $v$ in $\mathcal{U}_{ad},$ define $d(u,v)=P\otimes
dt\left\{  \left(  \omega,t\right)  \in\Omega\times\left[  0,T\right]  ;\text{
}u\text{(}\omega,t)\neq v\text{(}\omega,t)\right\}  ,$where $P\otimes dt$ is
the product of $P$ with the Lebesgue measure. It is clear that $d$ defines a
metric in $\mathcal{U}_{ad}.$

\begin{lemma}
\textit{i) (} $\mathcal{U}_{ad},d)$\textit{\ is a complete metric space.}

\textit{ii) For any} $p\geq1$ , there exists $M>0$ such that for any $u,v\in$
$\mathcal{U}_{ad}$,%
\[
E\left[  \sup\nolimits_{0\leq t\leq T}\left\vert X_{t}^{u}-X_{t}%
^{v}\right\vert ^{2p}\right]  \leq M.\left(  d(u,v)\right)  ^{1/2},
\]

\textit{where} $X_{t}^{u},$ $X_{t}^{v}$ \textit{are the solutions of
(}$\ref{MFSDE})$\textit{\ corresponding to} $u$ \textit{and} $v.$

\textit{iii) For any }$u$\textit{, }$v$\textit{\ in } $\mathcal{U}_{ad},$
$\left\vert J(u)-J(v)\right\vert \leq C.\left(  d(u,v)\right)  ^{1/2}.$
\end{lemma}

\bop The proof goes as in \cite{Zhou} Lemma 3.1 and uses classical arguments
from stochastic calculus, such as Burkholder-Davis-Gundy, H\"{o}lder's
inequality and Gronwall's lemma. The fact that the coefficients are of
mean-field type depending on the expectation of the state, does not add new
difficulties\eop

\noindent Let us define the Hamiltonian of the system associated to a random
variable $X:$

\begin{center}
$H(t,X,E(X),u,p,q)=b(t,X,E(X),u).p+\sigma(t,X,E(X),u).q-h(t,X,E(X),u)$
\end{center}

\noindent For any strict control $u\in\mathcal{U}$, we denote $(p,q)$ and
$(P,Q)$ the first and second order adjoint processes satisfying the following
backward SDEs%

\begin{equation}
\left\{
\begin{array}
[c]{ll}%
dp(t)= & -\left[  b_{x}(t)p(t)+E\left(  b_{y}(t)p(t)\right)  +\sigma
_{x}(t)q(t)+E\left(  \sigma_{y}(t)q(t)\right)  \right. \\
& \left.  -h_{x}(t)-E\left(  h_{y}(t)\right)  \right]  dt+q(t)dW_{t}+dM_{t}\\
p(T)= & -g_{x}(T)-E\left(  g_{y}(T)\right)
\end{array}
\right.  \label{Adjoint1}%
\end{equation}%
\begin{equation}
\left\{
\begin{array}
[c]{ll}%
-dP(t)= & -\left[  2b_{x}(t)P(t)+\sigma_{x}^{2}(t)P(t)+2\sigma_{x}%
(t)Q(t)+H_{xx}(t)]dt\right. \\
& +Q(t)dW_{t}+dN_{t}\\
P(T)= & -g_{xx}(x(T))
\end{array}
\right.  \label{Adjoint2}%
\end{equation}

\noindent where $X(t)$ is the state process associated with $u$,
$f_{x}(t)=f_{x}(t,X_{t},E(X_{t}),u_{t})$ for $f=b,$ $\sigma,$ $h$ and

\begin{center}
$H_{xx}(t,X,E(X),u,p,q)=b_{xx}(t,X,E(X),u).p+\sigma_{xx}(t,X,E(X),u).q-h_{xx}%
(t,X,E(X),u).$
\end{center}

\noindent$M$ and $N$ are square integrable martingales which are orthogonal to
the Brownian motion and are parts of the solutions. The appearance of such
martingales is due to the fact that $\left(  \mathcal{F}_{t}\right)  $ is not
necessarily the Brownian filtration. Note that BSDEs (\ref{Adjoint1}) and
(\ref{Adjoint2}) have been introduced for the first time in \cite{BucDjeLi},
without the orthogonal martingales $M$ and $N$.

\noindent Equation (\ref{Adjoint1}) is a mean field backward stochastic
differential equation (MFBSDE), whose driver is Lipschitz continuous, then by
\cite{BucLiPen} Theorem 3.1, it has a unique $\mathcal{F}_{t}-$adapted
solution $(p,q,M)$ such that:

\begin{center}
$E\left[  \sup\nolimits_{0\leq t\leq T}\left\vert p(t)\right\vert ^{2}+%
{\textstyle\int_{0}^{T}}
\left\vert q(t)\right\vert ^{2}dt+\left[  M,M\right]  _{T}\right]  <+\infty$
\end{center}

\noindent Note that in \cite{BucLiPen} Theorem 3.1, $\left(  \mathcal{F}%
_{t}\right)  $ is the Brownian filtration. Considering general filtrations on
which a Brownian motion is defined does not bring additional difficulties in
the proof of existence and uniqueness (see e.g \cite{ElkPenQuen} Theorem 5.1,
page 54).

\noindent Equation (\ref{Adjoint2}) is a\ classical backward stochastic
differential equation, whose driver is Lipschitz continuous, then by
\cite{ElkPenQuen} Theorem 5.1, it has a unique $\mathcal{F}_{t}-$adapted
solution $(P,Q,N)$ such that:

\begin{center}
$E\left[  \sup_{0\leq t\leq T}\left\vert P(t)\right\vert ^{2}+\int_{0}%
^{T}\left\vert Q(t)\right\vert ^{2}dt+\left[  N,N\right]  _{T}\right]
<+\infty.$
\end{center}

\noindent The following lemma is a stability result of the adjoints processes
with respect to the control variable.

\begin{lemma}
\textit{For any }$0<\alpha<1$\textit{\ and }$1<p<2$\textit{\ satisfying
}$(1+\alpha)<2,$\textit{\ there exists a constant }$C_{1}=C_{1}(\alpha
,p)>0$\textit{\ such that for any strict controls }$u$\textit{, }$u^{\prime}%
$\textit{\ along with the corresponding trajectories }$X$\textit{, }%
$X^{\prime}$\textit{\ and the solutions }$(p,q,P,Q,M,N)$\textit{, }%
$(p^{\prime},q^{\prime},P^{\prime},Q^{\prime},M^{\prime},N^{\prime}%
)$\textit{\ of the backward SDEs (\ref{Adjoint1}) and (\ref{Adjoint2}), the
following estimates hold}
\end{lemma}

\begin{center}
\textit{\ }$E\left[  \int_{0}^{T}\left(  \left\vert p(t)-p^{\prime
}(t)\right\vert ^{p}+\left\vert q(t)-q^{\prime}(t)\right\vert ^{p}\right)
dt+\left[  M-M^{\prime},M-M^{\prime}\right]  _{T}^{p/2}\right]  \leq
C_{1}d(u,u^{\prime})^{\frac{\alpha p}{2}}$

$E\left[  \int_{0}^{T}\left(  \left\vert P(t)-P^{\prime}(t)\right\vert
^{p}+\left\vert Q(t)-Q^{\prime}(t)\right\vert ^{p}\right)  dt+\left[
N-N^{\prime},N-N^{\prime}\right]  _{T}^{p/2}\right]  \mathit{\leq}%
C_{2}d(u,u^{\prime})^{\frac{\alpha p}{2}}$
\end{center}

\bop

\noindent The proof goes as in \cite{Zhou} Lemma 3.2. The only difference is
that the driver is of mean-field type. But this does not add new difficulties,
as the driver is linear and then Lipschitz in the state variable as well as in
its expectation.\eop

\noindent Let us define the $\mathcal{H-}$function or generalized Hamiltonian,
associated with a strict control $u$ and its state process $X,$ is defined as
follows:{}

\begin{center}
$\mathcal{H}^{(X(.),u(.))}(t,Y,E(Y),a)=H(t,Y,E\left(  Y\right)
,a,p(t),q(t)-P(t).\sigma(t,X_{t},E(X_{t}),u(t)))$

$-\frac{1}{2}\sigma^{2}(t,Y,E\left(  Y\right)  ,a)P(t)$
\end{center}

\noindent where $\left(  p(t),q(t)\right)  ,$ $\left(  P(t),Q(t)\right)  $ are
solutions of the adjoint equations (\ref{Adjoint1}) and (\ref{Adjoint2}).

\noindent According to the Chattering lemma and Proposition 3.9, for every
relaxed optimal control $\mu$, and for every $\varepsilon_{n}>0,$ there exist
a strict control $u^{n}$ such that:

\begin{center}
$J(u^{n})=J(\mu^{n})\leq\inf\left\{  J(\mu);\text{ }\mu\in\mathcal{R}\right\}
+\varepsilon_{n}.$
\end{center}

\noindent$u^{n}$ is called an $\varepsilon_{n}-$optimal\textit{ control. }

\noindent The next Proposition gives necessary conditions for near optimality
satisfied by an $\varepsilon_{n}-$optimal\textit{ control.}

\begin{proposition}
\textit{Let }$u^{n}$\textit{\ be an }$\varepsilon_{n}-$optimal\textit{ strict
control. Then there exist adapted processes }$(p^{n},q^{n},M^{n})$ and
$(P^{n},Q^{n},N^{n}),$ solutions of the adjoint equations (\ref{Adjoint1}) and
(\ref{Adjoint2}), corresponding to the admissible pair $(u^{n},X^{n})$ such
that\textit{:}%
\begin{equation}
E\left(
{\textstyle\int_{0}^{T}}
\mathcal{H}^{(X^{n}(t),u^{n}(t))}(t,X^{n}(t),E\left(  X^{n}(t)\right)
,u^{n}(t))dt\right)  \geq\text{ }\sup\nolimits_{a\in\mathbb{A}}E\left(
{\textstyle\int_{0}^{T}}
\mathcal{H}^{(X^{n}(t),u^{n}(t))}(t,X^{n}(t),E\left(  X^{n}(t)\right)
,a)dt\right)  -T\varepsilon_{n}^{1/3} \label{ApproxMaxPr}%
\end{equation}

\end{proposition}

\bop According to Lemma 4.2, the cost functional $J(u)$ is continuous with
respect to the topology induced by the metric $d.$ Then by applying Ekeland's
variational principle for $u^{n}$ with $\lambda_{n}=\varepsilon_{n}^{2/3}$,
there exists an admissible control $v^{n}$ such that $d(u^{n},v^{n}%
)\leq\varepsilon_{n}^{2/3}$ and $\widehat{J}(v^{n})\leq\widehat{J}(u)$ for all
$u\in\mathcal{U},$ where $\widehat{J}(u)=J(u)+\varepsilon_{n}^{1/3}%
d(u,v^{n}).$

The control $v_{n}$ which is $\varepsilon_{n}-$optimal is in fact optimal for
the new cost functional $\widehat{J}(u).$ We proceed as in the classical
mean-field maximum principle \cite{BucDjeLi} to derive a maximum principle for
$v^{n}.$ Let $t_{0}\in\left(  0,1\right)  ,$ $a\in\mathbb{A}$ and define the
spike variation of $v_{\delta}^{n}=a$ on $(t_{0},t_{0}+\delta)$ and
$v_{\delta}^{n}=$ $v_{n}(t)$ otherwise.

\noindent The fact that $\widehat{J}(v^{n})\leq\widehat{J}(u)$ for all
$u\in\mathcal{U}_{ad}$ and $d(v^{n},v_{\delta}^{n})\leq\delta$ imply that
$J(v_{\delta}^{n})-J(v^{n})\geq-\varepsilon_{n}^{1/3}\,\delta.$

\noindent Proceeding as in \cite{BucDjeLi}, we can expand $Y_{\delta}^{n}(.)$
(the solution of (\ref{MFSDE}) corresponding to $v_{\delta}^{n}$) to the
second order, to get the following inequality%

\[%
\begin{array}
[c]{l}%
E\left[  \int_{t_{0}}^{t_{0}+\delta}\frac{1}{2}\left(  \sigma(t,Y^{n}%
(t),E\left(  Y^{n}(t)\right)  ,a)-\sigma(t,Y^{n}(t),E\left(  Y^{n}(t)\right)
,v^{n}\right)  ^{2}P_{t}^{n}\right. \\
+p_{t}^{n}\left(  b(t,Y^{n}(t),E\left(  Y^{n}(t)\right)  ,a\right)
-b(t,Y^{n}(t),E\left(  Y^{n}(t)\right)  ,v^{n}))\\
+q_{t}^{n}\left(  \sigma(t,Y^{n}(t),E\left(  Y^{n}(t)\right)  ,a)-\sigma
(t,Y^{n}(t),E\left(  Y^{n}(t)\right)  ,v^{n})\right) \\
\left.  -\left(  h(t,Y^{n}(t),E\left(  Y^{n}(t)\right)  ,a\right)
-h(t,Y^{n}(t),E\left(  Y^{n}(t)\right)  ,v^{n}))dt\right]  +o(\delta
)\leq\varepsilon_{n}^{1/3}\,\delta
\end{array}
\]

\noindent where $Y^{n}(t)$ is the state process (solution of
(\ref{ApproxMaxPr})) corresponding to the control $v^{n}$ and ($p^{n},q^{n})$
and $(P^{n},Q^{n})$ are the first and second order adjoint processes,
solutions of (\ref{Adjoint1}) and (\ref{Adjoint2}) corresponding to
($v^{n},Y^{n})$.

\noindent The variational inequality is obtained for $v^{n}$ by dividing by
$\delta$ and tending $\delta$ to $0.$

\noindent The same claim can be proved for $u^{n}$ by using the stability of
the state equations and the adjoint processes with respect to the control
variable (Lemma 4.2\ and Lemma 4.3).\eop

\begin{remark}
The variational inequality (\ref{MFSDE}) can be proved with the supremum over
$a\in A$ replaced by the supremum over $u\in\mathcal{U}_{ad},$ by simply
putting $u(t)$ in place of $a$ in the definition of the strong perturbation.
\end{remark}

\subsection{The relaxed maximum principle}

We know that the relaxed control problem has an optimal solution $\mu.$
$\ $Let $X$ be the corresponding optimal state process. Let ($p,q,M)$ and
$(P,Q,N)$ the solutions of the first and second order adjoint equations,
associated with the optimal relaxed pair ($\mu,X)$.%

\begin{equation}
\left\{
\begin{array}
[c]{ll}%
dp(t)= & -\left[  \overline{b}_{x}(t)p(t)+E\left(  \overline{b}_{y}%
(t)p(t)\right)  +\overline{\sigma}_{x}(t)q(t)+E\left(  \overline{\sigma}%
_{y}(t)q(t)\right)  \right. \\
& \left.  -\overline{h}_{x}(t)-E\left(  \overline{h}_{y}(t)\right)  \right]
dt+q(t)dW_{t}+dM_{t}\\
p(T)= & -\overline{g}_{x}(T)-E\left(  \overline{g}_{y}(T)\right)
\end{array}
\right.  \label{RAdjoint1}%
\end{equation}%
\begin{equation}
\left\{
\begin{array}
[c]{ll}%
-dP(t)= & -\left[  2\overline{b}_{x}(t)P(t)+\overline{\sigma}_{x}%
^{2}(t)P(t)+2\overline{\sigma}_{x}(t)Q(t)+\overline{H}_{xx}(t)]dt\right. \\
& +Q(t)dW_{t}+dN_{t}\\
P(T)= & -\overline{g}_{xx}(x(T))
\end{array}
\right.  \label{RAdjoint2}%
\end{equation}

\noindent where we denote $\overline{f}(t)=f(t,x(t),E(x(t)),\mu(t))=\int
\nolimits_{A}f(t,x(t),E(x(t)),a)\mu(t,da)$ and $f$ stands for $b_{x}$,
$\sigma_{x},$ $h_{x},b_{y}$, $\sigma_{y},$ $h_{y},$ $H_{xx}.$

\noindent$\left(  M_{t}\right)  $ and $\left(  N_{t}\right)  $ are square
integrable martingales which are orthogonal to the Brownian motion $\left(
W_{t}\right)  .$

\noindent The drivers of the BSDEs (\ref{RAdjoint1}) and (\ref{RAdjoint2})
being Lipschitz continuous, then by \cite{BucLiPen} Theorem 3.1 and
\cite{ElkPenQuen} Theorem 5.1, they admit unique solutions ($p,q,M)$ and
$(P,Q,N)$ satisfying:

\begin{center}
$E\left[  \sup\nolimits_{0\leq t\leq T}\left\vert p(t)\right\vert ^{2}%
+\int_{0}^{T}\left\vert q(t)\right\vert ^{2}dt+\left[  M,M\right]
_{T}\right]  <+\infty,$

$E\left[  \sup\nolimits_{0\leq t\leq T}\left\vert P(t)\right\vert ^{2}%
+\int_{0}^{T}\left\vert Q(t)\right\vert ^{2}dt+\left[  N,N\right]
_{T}\right]  <+\infty.$
\end{center}

\noindent Define the generalized Hamiltonian function associated with the
optimal pair ($\mu,X(.))$ and the corresponding adjoint processes$,$

\begin{center}
$\mathcal{H}^{(X(.),\mu(.))}(t,Y,E(Y),a)=H(t,Y,E\left(  Y\right)
,a,p(t),q(t)-P(t).\sigma(t,X_{t},E(X_{t}),\mu(t)))$

$-\frac{1}{2}\sigma^{2}(t,Y,E\left(  Y\right)  ,a)P(t)$
\end{center}

\noindent The main result of this paper is the following.

\begin{theorem}
\textbf{(}Relaxed maximum principle\textbf{)}

\textit{Let }$(\mu,X)$\textit{\ be an optimal relaxed pair, then there exist
adapted processes }$(p,q,M)$\textit{ and }$(P,Q,N),$ solutions of the adjoint
equations (\ref{RAdjoint1}) and (\ref{RAdjoint2}) respectively, such that%
\begin{equation}
E\left(  \int_{0}^{T}\mathcal{H}^{(X(t),\mu(t))}(t,X(t),E(X(t)),\mu
(t))dt\right)  =\sup\nolimits_{a\in A}E\left(  \int_{0}^{T}\mathcal{H}%
^{(X(t),\mu(t))}(t,X(t),,E(X(t)),a)dt\right)  \label{PMax}%
\end{equation}

\end{theorem}

\noindent The proof of Theorem 4.6 will be given later.

\begin{corollary}
\textit{\ Under the same conditions as in Theorem 4.6 it holds that }%
\[
E\left(
{\textstyle\int_{0}^{T}}
\mathcal{H}^{(X(t),\mu(t))}(t,X(t),E(X(t)),\mu(t))dt\right)  =\sup
\nolimits_{\upsilon\in\mathbb{P}\left(  A\right)  }E\left[
{\textstyle\int_{0}^{T}}
\mathcal{H}^{(X(t),\mu(t))}(t,X(t),E(X(t)),\upsilon)\right]  dt
\]

\end{corollary}

where $\mathcal{H}^{(X(t),\mu(t))}(t,X(t),,E(X(t)),\upsilon)=\int_{\mathbb{A}%
}\mathcal{H}^{(X(t),\mu(t))}(t,X(t),,E(X(t)),a)\upsilon(da)$ and
$\mathbb{P}\left(  \mathbb{A}\right)  $ is the space of probability measures
on $\mathbb{A}$.

\bop The inequality

\begin{center}
$\sup\nolimits_{\mu\in\mathbb{P}\left(  \mathbb{A}\right)  }E(\int_{0}%
^{T}\mathcal{H}^{(X(t),\mu(t))}(t,X(t),E(X(t)),\upsilon)dt)\geq\sup
\nolimits_{a\in A}E\left(  \int_{0}^{T}\mathcal{H}^{(X(t),\mu(t))}%
(t,X(t),E(X(t)),a)dt\right)  $
\end{center}

\noindent is obvious. Indeed it suffices to take $\mu=\delta_{a}$ where $a$ is
any element of $\mathbb{A}.$ Now if $\upsilon\in\mathbb{P}\left(  A\right)  $
is a probability measure on $\mathbb{A},$ then

\begin{center}
$\int_{0}^{T}\mathcal{H}^{(X(t),\mu(t))}(t,X(t),E(X(t)),\upsilon)dt\in
conv\left\{  E\left(  \int_{0}^{T}\mathcal{H}^{(X(t),\mu(t))}%
(t,X(t),,E(X(t)),a)dt\right)  ,a\in\mathbb{A}\right\}  $
\end{center}

\noindent Hence, by using Fubini's theorem, it holds that

\begin{center}
$\int_{0}^{T}\mathcal{H}^{(X(t),\mu(t))}(t,X(t),E(X(t)),\upsilon)dt\leq
\sup\nolimits_{a\in A}E\left(  \int_{0}^{T}\mathcal{H}^{(X(t),\mu
(t))}(t,X(t),,E(X(t)),a)dt\right)  $
\end{center}

\eop

\textbf{Remark .} Since $\mathbb{P}(\mathbb{A})$ is a subspace of $\mathbb{V}$
whose elements are constant (in $\left(  \omega,t\right)  )$ relaxed controls,
then one can replace in Corollary 4.8, the supremum over $\mathbb{P}\left(
A\right)  $ by the supremum over $\mathbb{V}.$

\begin{corollary}
(The Pontriagin \textit{relaxed maximum principle). If }$\left(  \mu,X\right)
$ \textit{denotes an optimal relaxed pair, then there exists a Lebesgue
negligible subset }$N$\textit{\ such that for any }$t$\textit{ not in }%
$N$\textit{\ \ }%
\[
\mathcal{H}^{(X(t),\mu(t))}(t,X(t),E(X(t)),\mu(t))=\sup\nolimits_{\mu
\in\mathbb{P}\left(  \mathbb{A}\right)  }\mathcal{H}^{(X(t),\mu(t))}%
(t,X(t),E(X(t)),\upsilon),\text{ }P-a.s.
\]

\end{corollary}

\bop Let $\theta\in\left]  0,T\right[  $ and $B\in\mathcal{F}_{\theta},$ for
small $h>0$ define the relaxed control $\mu_{t}^{h}=\upsilon1_{B}$ \ \ for
$\theta<t<\theta+h$ and $\mu_{t}^{h}=\mu_{t}$ \ \ otherwise, where $\upsilon$
is a probability measure on $A.$ It follows from \textit{Theorem 4.6,} that

\begin{center}
$1/h\int_{\theta}^{\theta+h}\mathit{E}\left[  1_{B}\mathcal{H}^{(X(t),\mu
(t))}(t,X(t),E(X(t)),\mu(t))\right]  dt\geq\;1/h\int_{\theta}^{\theta
+h}\mathit{E}\left[  1_{B}\mathcal{H}^{(X(t),\mu(t))}(t,X(t),E(X(t)),\upsilon
)\right]  dt$
\end{center}

\noindent Therefore passing to the limit as $h$ tends to zero, we obtain

\begin{center}
$\mathit{E}\left[  1_{B}\mathcal{H}^{(X(\theta),\mu(\theta))}(\theta
,X(\theta),E(X(\theta)),\mu(\theta))\right]  \geq\mathit{E}\left[
1_{B}\mathcal{H}^{(X(\theta),\mu(\theta))}(\theta,X(\theta),E(X(\theta
)),\upsilon)\right]  $
\end{center}

\noindent for any $\theta$ not in some Lebesgue null set $N.$

\noindent The last inequality is valid for all $B\in\mathcal{F}_{\theta},$
then for any bounded $\mathcal{F}_{\theta}$-measurable random variable $F,$ we get

\begin{center}
$E\left[  F\mathcal{H}^{(X(t),\mu(t))}(t,X(t),E(X(t)),\mu(t))\right]  \geq
E\left[  F\mathcal{H}^{(X(t),\mu(t))}(t,X(t),E(X(t)),\upsilon)\right]  $
\end{center}

\noindent which leads to

\begin{center}
$E\left[  \mathcal{H}^{(X(\theta),\mu(\theta))}(\theta,X(\theta),E(X(\theta
)),\mu(\theta))/\mathcal{F}_{\theta}\right]  \geq E\left[  \mathcal{H}%
^{(X(\theta),\mu(\theta))}(\theta,X(\theta),E(X(\theta)),\upsilon
)/\mathcal{F}_{\theta}\right]  $
\end{center}

\noindent The result follows from the measurability with respect to
$\mathcal{F}_{\theta}$ of the quantities inside the conditional
expectation.\eop

\noindent The proof of Theorem 4.6 is based on the following lemma.

\begin{lemma}
\textit{Let }$\left(  p^{n},q^{n}\right)  $, $(P^{n},Q^{n})$\textit{\ (resp.
}$\left(  p,q\right)  $, $(P,Q)$) \textit{be the solutions of the first and
second order adjoint equations (\ref{Adjoint1}) and (\ref{Adjoint2})}
associated with the pair $(u^{n},X^{n}),$ (resp. solutions of first and second
order adjoint equations (\textit{\ref{RAdjoint1}) and (\ref{RAdjoint2})
associated to (}$\mu,X)$)\textit{, then }
\end{lemma}

\begin{center}
\textit{i) }$\lim\limits_{n\rightarrow+\infty}E\left[  \int_{0}^{T}\left(
\left\vert p^{n}(t)-p(t)\right\vert ^{2}+\left\vert q^{n}(t)-q(t)\right\vert
^{2}\right)  dt+\left[  M^{n}-M,M^{n}-M\right]  _{T}\right]  =0$

\textit{ii) }$\lim\limits_{n\rightarrow+\infty}E\left[  \int_{0}^{T}\left(
\left\vert P^{n}(t)-P(t)\right\vert ^{2}+\left\vert Q^{n}(t)-Q(t)\right\vert
^{2}\right)  dt+\left[  N^{n}-N,N^{n}-N\right]  _{T}\right]  =0$

\textit{iii) }$\lim\limits_{n\rightarrow+\infty}E\left(  \int_{0}%
^{T}\mathcal{H}^{(X^{n}(t),u^{n}(t))}(t,X^{n}(t),E\left(  X^{n}(t)\right)
,u^{n}(t))dt\right)  =E\left(  \int_{0}^{T}\mathcal{H}^{(X(t),\mu
(t))}(t,X(t),E\left(  X(t)\right)  ,\mu(t))dt\right)  $
\end{center}

\bop i)Let us write down the drivers of the first order adjoint equations
(\textit{\ref{Adjoint1}) and (}\ref{RAdjoint1})\textit{ }corresponding to
$(u^{n},X^{n})$ and ($\mu,X).$

\begin{center}
$G^{n}(t,p_{t}^{n},q_{t}^{n})=-b_{x}^{n}(t)p^{n}(t)+E\left(  b_{y}^{n}%
(t)p^{n}(t)\right)  +\sigma_{x}^{n}(t)q^{n}(t)+E\left(  \sigma_{y}^{n}%
(t)q^{n}(t)\right)  -h_{x}^{n}(t)-E\left(  h_{y}^{n}(t)\right)  $

$G(t,p_{t},q_{t})=-\overline{b}_{x}(t)p(t)+E\left(  \overline{b}%
_{y}(t)p(t)\right)  +\overline{\sigma}_{x}(t)q(t)+E\left(  \overline{\sigma
}_{y}(t)q(t)\right)  -\overline{h}_{x}(t)-E\left(  \overline{h}_{y}(t)\right)
$
\end{center}

\noindent where

$f^{n}(t)=f(t,X_{t}^{n},E(X_{t}^{n}),u_{t}^{n})=%
{\textstyle\int_{\mathbb{A}}}
f(t,X_{t}^{n},E(X_{t}^{n}),a)\delta_{u_{t}^{n}}(da)$ for $f=b_{x},$
$\sigma_{x},$ $h_{x},b_{y},\sigma_{y},h_{y}.$

$\overline{f}(t)=f(t,X(t),E(X_{t}),\mu(t))=%
{\textstyle\int_{\mathbb{A}}}
f(t,X(t),E(X_{t}),a)\mu(t,da)$ where $f$ stands for $b_{x}$, $\sigma_{x},$
$h_{x},$ $b_{y},\sigma_{y},h_{y}.$

\noindent By using a subtil stability result of Hu and Peng \cite{HuPeng},
Theorem 2.1, it is sufficient to show that:

\begin{center}
$\underset{n\rightarrow\infty}{\lim}E\left[  \left\vert \int_{t}^{T}\left(
G^{n}(t,p_{t},q_{t})-G(t,p_{t},q_{t})\right)  dt\right\vert ^{2}\right]  =0$
\end{center}

\noindent We have

\begin{center}%
\begin{equation}%
\begin{array}
[c]{cl}%
\left\vert \int_{t}^{T}\left(  G^{n}(t,p_{t},q_{t})-G(t,p_{t},q_{t})\right)
dt\right\vert  & \leq\left\vert \int_{t}^{T}\left(  b_{x}^{n}(t)-\overline
{b}_{x}(t)\right)  p(t)dt\right\vert +\left\vert \int_{t}^{T}E\left[  \left(
b_{y}^{n}(t)-\overline{b}_{y}(t)\right)  p(t)\right]  dt\right\vert \\
& +\left\vert \int_{t}^{T}\left(  \sigma_{x}^{n}(t)-\overline{\sigma}%
_{x}(t)\right)  q(t)dt\right\vert +\left\vert \int_{t}^{T}E\left[  \left(
\sigma_{y}^{n}(t)-\overline{\sigma}_{y}(t)\right)  q(t)\right]  dt\right\vert
\\
& +\left\vert \int_{t}^{T}\left(  h_{x}^{n}(t)-\overline{h}_{x}(t)\right)
dt\right\vert +\left\vert \int_{t}^{T}E\left(  h_{y}^{n}(t)-\overline{h}%
_{y}(t)\right)  dt\right\vert
\end{array}
\label{InequalityDrivers}%
\end{equation}

\end{center}

\noindent Let us treat the first term in the right hand side of
(\ref{InequalityDrivers}).

\begin{center}
$%
\begin{array}
[c]{cc}%
\int_{t}^{T}\left(  b_{x}^{n}(t)-\overline{b}_{x}(t)\right)  p(t)dt &
=\int_{t}^{T}\left(
{\textstyle\int_{\mathbb{A}}}
b_{x}(t,X_{t}^{n},E(X_{t}^{n}),a)\delta_{u_{t}^{n}}(da)-%
{\textstyle\int_{\mathbb{A}}}
b_{x}(t,X_{t},E(X_{t}),a)\mu_{t}(da)\right)  p(t)dt\\
& =\int_{t}^{T}\left(
{\textstyle\int_{\mathbb{A}}}
b_{x}(t,X_{t}^{n},E(X_{t}^{n}),a)\delta_{u_{t}^{n}}(da)-%
{\textstyle\int_{\mathbb{A}}}
b_{x}(t,X_{t},E(X_{t}),a)\delta_{u_{t}^{n}}(da)\right)  p(t)dt\\
& +\int_{t}^{T}\left(
{\textstyle\int_{\mathbb{A}}}
b_{x}(t,X_{t},E(X_{t}),a)\delta_{u_{t}^{n}}(da)-%
{\textstyle\int_{\mathbb{A}}}
b_{x}(t,X_{t},E(X_{t}),a)\mu_{t}(da)\right)  p(t)dt
\end{array}
$
\end{center}

\noindent The facts that $b_{x}$ is Lipschitz continuous in $(x,y)$ and
$\left(  X_{t}^{n}\right)  $ converges to $X_{t}$ uniformly in $t$ in
probability, imply that the first term in the right hand side of the last
inequality converges to $0$ in probability$.$

\noindent Furthermore $E\left(  \sup_{0\leq t\leq T}\left\vert p(t)\right\vert
^{2}\right)  <+\infty,$ then $\sup_{0\leq t\leq T}\left\vert p(t)\right\vert
<+\infty,$ $P-a.s,$ which implies the existence of a $P-$negligible set $N,$
such that for each $\omega\notin N,$ there exist $M(\omega)<+\infty$ s.t
$\sup_{0\leq t\leq T}\left\vert p(t)\right\vert \leq M(\omega).$ In particular
for each $\omega\notin N,$ the function $b_{x}(t,X_{t},E(X_{t}%
),a)p(t).1_{\left[  0,t\right]  }$ is a measurable bounded function in $(t,a)$
and continuous in $a$, therefore it is a test function for the stable
convergence$.$ Hence by using the fact that $\left(  \delta_{u_{t}^{n}%
}(da)\,dt\right)  $\textit{\ converges in }$\mathbb{V}$\textit{ to }$\mu
_{t}(da)\,dt,$\textit{\ }$P-a.s.$, it follows that the second term in the
right hand side tends to $0,$ $P-a.s$.

\noindent We conclude by using the Lebesgue dominated convergence theorem.

\noindent The other terms containing $p(t)$ can be handled by using the same techniques.

\noindent The terms in (\ref{InequalityDrivers}) containing $q(t)$ can be
treated similarily. However one should pay a little more attention as $q(t)$
is just square integrable (in $(t,\omega)$). More precisely

\begin{center}
$%
\begin{array}
[c]{cc}%
\left\vert \int_{t}^{T}\left(  \sigma_{x}^{n}(t)-\overline{\sigma}%
_{x}(t)\right)  q(t)dt\right\vert  & \leq\left\vert \int_{t}^{T}\left(
\sigma_{x}^{n}(t)-\overline{\sigma}_{x}(t)\right)  q(t)1_{\left\{  \left\vert
q(t)\right\vert \leq N\right\}  }dt\right\vert +\left\vert \int_{t}^{T}\left(
\sigma_{x}^{n}(t)-\overline{\sigma}_{x}(t)\right)  q(t)1_{\left\{  \left\vert
q(t)\right\vert \geq N\right\}  }dt\right\vert
\end{array}
$
\end{center}

\noindent The first integral in the right hand side may be handled by using
similar argument as precedently as the function $\left(  \sigma_{x}%
^{n}(t)-\overline{\sigma}_{x}(t)\right)  q(t)1_{\left\{  \left\vert
q(t)\right\vert \leq N\right\}  }$ is measurable bounded and continuous in
$a.$ The second term tends to $0$ by Tchebychev's inequality, using the square
integrability of $q(t).$

\noindent ii) and iii) are proved by using similar arguments.\eop

\bop\textbf{of Theorem 4.9. }The result is proved by passing to the limit in
inequality (\ref{ApproxMaxPr}) and using Lemma 4.9.\eop

\textbf{Acknowledgements.} The authors would like to thank the referees for
many helpful suggestions, which lead to a substancial improvement of the manuscript.

\end{document}